\documentclass[12pt]{article}
\usepackage{amssymb,amsmath,cite,bm}
\usepackage{psfrag,subfigure}
\newtheorem{remark}{Remark}

\newtheorem{properties}{Properties}

\begin{document}

\title{\bf Nonsmooth Mechanics Based on Linear Projection Operator}

\author{Farhad Aghili\thanks{email:faghili@encs.concordia.ca}}

\date{}
\maketitle

\begin{abstract}
This paper presents a unifying dynamics formulation for nonsmooth multibody systems (MBSs) subject to changing topology and multiple impacts based on linear projection operator. An oblique projection matrix  ubiquitously derives all characteristic variables of such systems  as follow: i) The constrained acceleration before jump discontinuity  from projection of unconstrained acceleration, ii) post-impact velocity from projection of pre-impact velocity, iii) impulse  during impact from projection of pre-impact momentum, iv) generalized constraint force from projection of  generalized input force, and v) post-impact kinetic energy from  pre-impact kinetic energy based on projected inertia matrix. All solutions are presented in closed-form with elegant geometrical interpretations. The formulation is general enough to be applicable to MBSs subject to simultaneous multiple impacts with non-identical restitution coefficients, changing topology, i.e., unilateral constraints becomes inactive or vice versa, or even when the overall constraint Jacobian becomes singular. Not only do the solutions always exist regardless the constraint condition, but  also the condition number for a generalized constraint inertia matrix  is minimized in order to reduce numerical sensitivity in computation of the projection matrix to roundoff errors. The model is proven to be energetically consistent if a global restitution coefficient is assumed.  In the case of non-identical restitution coefficients, the set of energetically consistent restitution matrices is characterized by using Linear Matrix Inequality (LMI).
\end{abstract}

\section{Introduction}

The class of robotic systems with varying topology due to switching constraints or impact event arises in many robotics applications involving multiple contacts or formations  of closed-loop topology \cite{McClamorch-Wang-1988,GarciadeJalon-Bayo-1994,Aghili-Buehler-Hollerbach-1997a,Blajer-Schiehlen-Schirm-1994,Aghili-2010h,Aghili-Su-2017} such as  robotic assembly \cite{Gottschlich-Kak-1989}, force control of constrained robots  \cite{Dupree-Liang-2008,Aghili-2020a}, and legged robotics \cite{Marhefka-Orin-1999}. These robotic systems can  be generally treated as multibody systems (MBSs) subject to unilateral kinematic constraints with a time-varying topology. In many cases, the constraints arises from intermittent contacts, which can fundamentally change the dynamic behavior of the system e.g., when passive unilateral constraints become active and vice versa. There are two main challenges in dynamics formulation of these systems: i) time-varying structure of the systems due to activation and deactivation of unilateral constraints, ii)  non-smooth  behaviors associated with the jump discontinuity in the velocities  and impulsive constraint force during constraint activation events. More specifically, the unilateral constraints are defined so that they are active when the corresponding relative displacements  and constraint forces are non-negative and otherwise they remain inactive. Consequently the dimension of the overall constraint manifold, which is equal to the  number of independent kinematic constraints, varies over time and  so does the number of  degrees-of-freedom exhibited by the system. Moreover when one or more inactive unilateral constraints become active, such events transpire impact, which gives rise to impulsive constraint forces. The topology change is often modeled as instantaneous events meaning that  the system's velocities cannot be assumed continuous  during constraint switching or impact events. In that case, the system states can not be calculated from an acceleration model, instead they have to be derived from an impact model by incorporating the impulse-momentum balance of the entire system together with a restitution law. Simulation and analysis of such constrained robotic systems with varying topology  call for a unifying formulation that be used not only for a single motion dynamics model of MBSs with varying topology but for the impulsive motion problem as well.  Ideally, a unified formulation which yields closed-form solution is desirable  not only for computational efficiency but for giving geometrical interpretation of multiple impact phenomenon in MBSs.

In the literature, the primary approach to modeling impact in MBSs is based on combining the equation of momentum-balance together with a restitution law \cite{Khulief-2013,Brogliato-2014}. In this approach, the impact is assumed to be an infinitesimal event and subsequently the set of the impulse-momentum equations are algebraically solved to produce velocity jump rather than updating the velocity from  integration of the acceleration vector. In this formulation, the impact is characterized by the coefficient of restitution, which is defined as the ratio of the local velocities after and before collision.   Other approaches assume smooth compliant modelling of impact where  the impulsive force is typically presented by a linear or nonlinear spring-damper model \cite{Khulief-2013,Hunt-Crossley-1975}, e.g., the Haunt-Crossley nonlinear spring-dashpot model. However, it has been established that   a linear or non-linear  compliant model becomes equivalent to the momentum-balance model if the coefficient of restitution is specifically selected according to the damping and stiffness properties of the compliance model \cite{Goldsmith-1960,Lankarani-Nikravesh-1990,Marhefka-Orin-1999}. Impact model for  collisions of a special class of planar kinematic chains is presented in \cite{Hurmuzlu-Chang-1992}. Impact dynamics of a five-link biped walking on level ground are studied in \cite{Mu-Wu-2006}  for correlating the gait parameters with the contact event following impact. Implementation and validation of theoretical impact model representing non-smooth dynamics  exhibited in a bouncing dimer in a planar environment are reported in \cite{Liu-Zhao-2008,Zhao-Liu-2009}. Impact models for humanoid robots for respective controllers are developed in \cite{Yoshida-Takeuchi-2014} to ensure smooth reaction control under both impact-force and continuous-force disturbances. Although it is common practice to treat the coefficient of restitution  to be independent of the impact velocity, empirical observations have demonstrated otherwise  \cite{Hunt-Crossley-1975,Schiehlen-Seifried-2006}. Dependency of the coefficient of restitution to the local approach velocity implies that  a global coefficient of restitution may not be assumed for concurrent simulation of for multiple impacts even with identical material properties of the contact surfaces \cite{Lankarani-Nikravesh-1990,Marhefka-Orin-1999,Zhang-Vu-Quoc-2002,Najafabadi-Kovecses-Angeles-2008,Aghili-2019a}.

Another challenge in modelling mechanical systems subject to unilateral constraints is that they often have a time-varying topology. This is because when passive unilateral constraints becomes active the dynamics behavior of the system  changes and so do their numbers of degrees-of-freedom. This means that it is not possible to find beforehand a global minimal-order model characterized with an  independent coordinates having a fixed dimension. Utilizing a set of non-minimal models corresponding to different constraint situation necessitates switching between different models during simulation or control and thus this approach suffers from  smoothness and stability of the transition \cite{McClamorch-Wang-1988,GarciadeJalon-Bayo-1994,Aghili-Buehler-Hollerbach-1997a,Blajer-Schiehlen-Schirm-1994}.
Another shortcoming of reduced-order dynamics model is  the non-unique  relationship between the independent coordinates and the spatial configuration of MBSs. Moreover, using the minimum-order independent coordinates becomes particularly limiting for MBSs  passing through singular configurations, which inevitably  gives rise to the numbers of DOF. Regularization of singularities  of bilaterally constrained systems is studied in the literature \cite{Arponen-2001,Blajer-2002b,Aghili-Piedboeuf-2003a,Muller-2006}. Explicit dynamics formulation in terms of dependent coordinates was proposed in our earlier works  based on the notion of orthogonal projection that does not require  a recourse to an iterative process \cite{Aghili-Piedboeuf-2003a,Aghili-2005}. Subsequently, similar projection approaches has been used for modelling and control of robots with switching topology such as walking robots  \cite{Mistry-Buchli-Schaal-2010,Mistry-Righetti-2011,Righetti-Buchli-Mistry-2011,Aghili-2019a}.

This paper presents a unifying approach for modelling nonsmooth mechanics based on oblique projection operator. More specifically, the acceleration of constrained MBS is simply derived by  projecting the acceleration computed from the same MBS without imposing the constraints. Likewise, the velocity jump and impulse during impact are directly obtained by  obliquely projecting the pre-impact velocity and pre-impact momentum of the MBS, respectively. Other post-impact variables such as impact and kinetic energy are elegantly derived based on projection.  Therefore, the dynamics model seemingly works not only when  the unilateral constraints remains active through a finite time interval but when changing topology occurs, i.e., unilateral constraints becomes inactive or vice versa, or even when the overall constraint Jacobian becomes singular. The other advantages of the formulation are it can handle  MBSs with chain tree and close-loop topologies or subject to  simultaneous multiple contacts  with non-identical restitution coefficients \cite{Aghili-2019a}. The formulation is energetically consistent if a global restitution coefficient can be assumed, while  for the case of non-identical restitution coefficients the energetically consistency is satisfied upon a LMI condition.

\section{Mechanical System with Bilataeral/Unilateral Constraints}
Consider a mechanical system  with generalized coordinate $\bm q \in \mathbb{R}^n$ subject to a set of $m_b$ bilateral constrains and $m_u$ unilateral constraints. Dynamics equation of such system can be described by
\begin{subequations} \label{eq:dynamics_unilateral}
\begin{align} \label{eq:Mddotq}
\bm M(\bm q) \ddot{\bm q}  + \bm h(\bm q, \dot{\bm q})  = & \bm u + \bm A_b^
T \bm\lambda_b +  \bm A_u^T \bm\lambda_u\\ \label{eq:phib}
\bm\phi_b(\bm q) &= \bm 0 \\ \label{eq:phiu}
\bm\phi_u(\bm q) &\geq \bm 0, \qquad \bm\lambda_u \geq \bm 0 \\ \label{eq:lam>0}
\bm\phi_u^T(\bm q) \bm\lambda_u & = 0
\end{align}
\end{subequations}
where ${\bm  M}(\bm q) \in \mathbb{R}^{n \times n}$  is the inertia matrix;
${\bm h}({\bm q}, \dot{\bm  q}) \in\mathbb{R}^{n}$ contains Coriolis,
centrifugal, and gravitational terms, and $\bm u \in\mathbb{R}^{n}$ represents collectively all actuation, dissipative, or external  forces acting on the system.  In the above equation, $\bm\phi_b(\bm q)\in \mathbb{R}^{m_b}$ and $\bm\phi_u(\bm q)\in \mathbb{R}^{m_u}$ are the set of bilateral and unilateral constraint equations;  ${\bm\lambda}_b \in \mathbb{R}^{m_b}$ and ${\bm\lambda}_u \in \mathbb{R}^{m_u}$ are the Lagrangian multipliers  associated with the bilateral and unilateral constraints, and \eqref{eq:lam>0} is the complementarity condition. The jacobians associated with the constraints are
\begin{equation} \label{eq:AbAu}
\bm A_b = \frac{\partial \bm\phi_b}{\partial \bm q} \quad \mbox{and} \quad
\bm A_u = \bm\Gamma \frac{\partial \bm\phi_u}{\partial \bm q}
\end{equation}
where $\bm\Gamma$ is a diagonal activation matrix whose
entries are computed according to
\begin{equation}
\Gamma_{ii} =\left\{ \begin{array}{ll} 1 & \mbox{if} \quad \phi_{u_i} = 0 \;  \wedge \; \lambda_{u_i} \geq 0
\\ 0 & \mbox{otherwise} \end{array} \right.
\end{equation}
In other words, $\Gamma_{ii}=1$ implies activation of the $i$th unilateral constraint. Notice that \eqref{eq:AbAu} simply ignores un-activated unilateral constraints as if they did not exist. This is because  if a unilateral constraint is not activated, then the corresponding switch $\Gamma_{ii}$ is zero making the constraint passive. As a result, the rank of Jacobian matrix $\bm A_u$ increases or decreases as some of the passive unilateral constraints become active or vice versa. As shall be described in  Section~\ref{sec:impact}, during
the transition when the rank abruptly  changes,  the continuity may not be ensured and therefore the system should be cast in the framework of non-smooth mechanics.

Denote the entire vector of Lagrangian multiplier by $\bm\lambda=[\bm\lambda_b^T \;\; \bm\lambda_u^T ]^T$ and the augmented Jacobian matrix
\begin{equation} \label{eq:Adef}
\bm A = \begin{bmatrix} \bm A_b \\ \bm A_u \end{bmatrix}, \quad \mbox{and hence} \quad \bm f = \bm A^T \bm\lambda
\end{equation}
is the generalized constraint force. The  equations of the constrained mechanical system
can be transcribed as in the standard dynamics equations having only equality constraints
\begin{subequations}
\begin{equation} \label{eq:Mddotq_standard}
\bm M(\bm q) \ddot{\bm q}  + \bm h(\bm q, \dot{\bm q})=  \bm u + \bm f,
\end{equation}
\begin{equation}\label{eq:Adq}
\bm A \dot {\bm q}  =\bm 0
\end{equation}
\end{subequations}
The Jacobian matrix may not be full-rank, that is $\text{rank}(\bm A)=r$ and
$r\leq m$, where $m=m_u+m_b$ is the sum of all constraints. Alternatively, the constraint condition \eqref{eq:Adq} can be equivalently described by  $\bm P \dot{\bm q}= \dot{\bm q}$ \cite{Aghili-Piedboeuf-2003a}, where the square matrix $\bm P \in \mathbb{R}^{n \times n}$ represents the {\em orthogonal projection} onto the null-space of $\bm A$, and hence $\bm P^2= \bm P = \bm P^T$. The projection matrix can be obtained from $\bm P \triangleq \bm I -  \bm A^+ \bm A$ where $\bm A^+$ is the pseudo-inverse of $\bm A$. Hereafter, ${\cal N}$ and ${\cal N}^{\perp}$ represent, respectively, the null-space and null-space complement of $\bm A$.
As shown in Appendix~\ref{appx_skew}, one can also show that the generalized velocity is linearly mapped to  the ${\cal N}^{\perp}$ component of the generalized acceleration $\ddot{\bm q}_{\perp} := ({\bm I - \bm P}) \ddot{\bm q}$ through
\begin{equation}\label{eq:Omega_dq}
\ddot{\bm q}_{\perp} = \bm\Omega \dot{\bm  q},
\end{equation}
where $\bm\Omega =\bm\Lambda^T - \bm\Lambda$ is a skew-symmetric, and $\bm\Lambda = - \bm A^+ \dot{\bm  A} \bm P$.

On the other hand, since  $\bm P \bm f = \bm 0$, the constraint force can be simply eliminated from
equation \eqref{eq:Mddotq_standard} if both sides of the latter equation are
pre-multiplied by $\bm P$, i.e.,
\begin{equation} \notag
\bm P \bm M \ddot{\bm q}  = \bm P (\bm u -\bm h),  \qquad \mbox{or}
\end{equation}
\begin{equation} \label{eq:PMPddq}
\bm P \bm M \bm P \ddot{\bm q} = \bm P (\bm u - \bm h)  -  \bm P \bm M \bm\Omega \dot{\bm q},
\end{equation}
where $\bm h$ denotes $\bm h(\bm q, \dot{\bm q})$ for brevity of notation. Pre-multiplying equation \eqref{eq:Omega_dq} by positive scalar $\nu>0$, and then add both sides of the resulting equation with \eqref{eq:PMPddq}, we arrive at
\begin{equation} \label{eq_EqMotion2}
{\bm M}_c \ddot{\bm q} = \big( \nu \bm I - \bm P \bm M \big) \bm\Omega \dot{\bm q}  +  \bm P (\bm u - \bm h),
\end{equation}
where ${\bm M}_c =   \bm P  \bm M  \bm P + \nu (\bm I - \bm P)$ is the constraint inertia matrix. Equation \eqref{eq_EqMotion2} constitutes the {\em non-minimal order} dynamics model of mechanical systems in the standard form where $\bm M_c(\bm q)$ can be treated as the constraint inertia matrix.
\begin{properties} \label{prop:Mc}
The constraint inertia matrix satisfies the following properties:
\begin{subequations}
\begin{eqnarray} \label{eq:Mc>0}
& \bm M_c \succeq 0  \\  \label{eq:PMP}
& \bm P \bm M_c = \bm M_c \bm P = \bm P \bm M \bm P \\ \label{eq:Mc_commute}
& \bm P \bm M_c^{-1} = \bm M_c^{-1} \bm P = \bm P \bm M_c^{-1} \bm P = \bm M_o^+\\ \label{eq:MMcinv}
& \bm M \bm M_c^{-1} = \bm P \\ \label{eq:McinvOmega}
& \nu \bm M_c^{-1} \bm\Omega \dot{\bm q} =  \bm\Omega \dot{\bm q}\\ \label{eq:condMc}
& \min_{\nu} {\rm cond}(\bm M_c) \quad  \Leftarrow \quad   \lambda_{\stackrel{\rm min}{\neq 0}}(\bm M_o) \leq \nu \leq \lambda_{\rm max}(\bm M_o)
\end{eqnarray}
\end{subequations}
where $\bm M_o=\bm P \bm M \bm P$.
\end{properties}

That $\bm P$ commutes with both the constraint inertia matrix and its inverse, i.e., \eqref{eq:Mc_commute}, can be inferred readily from definition of the constraint inertia matrix. The proofs of that matrix $\bm M_c$ is always  symmetric positive-definite \eqref{eq:Mc>0}, identity \eqref{eq:McinvOmega}, and minimizing the condition number of the constraint inertia matrix  are given in Appendix~\ref{appx:Mc_properties}. Notice that condition \eqref{eq:condMc} is important from a numerical point of view as computation of a projection matrix defined in the following section  will require inversion of the constraint inertia matrix.

\subsection{Generalized acceleration through oblique projection}
Let us define
\begin{equation}  \label{eq:S_definition}
\bm S :=\bm I- \bm M_c^{-1}\bm P \bm M=\bm I- \bm M_o^+ \bm M,
\end{equation}
which is not generally a symmetric matrix, i.e., $\bm S^T \neq \bm S$. It can be inferred from properties \eqref{eq:Mc_commute} and \eqref{eq:MMcinv} that $\bm P \bm M \bm P\bm M_c^{-1} \bm P = \bm P \bm M \bm M_c^{-1} \bm P = \bm P$. Using the latter identity in the following derivation leads to
\begin{equation} \notag
\bm S^2 = \bm I -2 \bm M \bm M_c^{-1} \bm P + \bm M \bm M_c^{-1} \underbrace{\bm P \bm M \bm P \bm M_c^{-1} \bm P}_{\bm P} = \bm S.
\end{equation}
Therefore, since $\bm S^2 = \bm S$ and $\bm S$ is not a symmetric matrix then $\bm S$ must be an {\em oblique projector}.
\begin{properties}
It can be verified that projection matrix $\bm S$ satisfies the following identities
\begin{subequations} \label{eq:S_prop}
\begin{eqnarray} \label{eq:SP=0}
& \bm S \bm P = \bm P \bm S^T = \bm 0 \\ \label{(I-P)S}
& \bm S (\bm I - \bm P)= \bm S \\ \label{(I-P)St}
& (\bm I - \bm P) \bm S= \bm I - \bm P \\ \label{MS=StM}
& \bm M \bm S  = \bm S^T \bm M = \bm S^T \bm M \bm S\\ \label{eq:(I-S)invM}
& (\bm I - \bm S) \bm M^{-1} =\bm M^{-1}(\bm I - \bm S^T) = \bm M_o^+
\end{eqnarray}
\end{subequations}
\end{properties}
See Appendix~\ref{appx:S_properties} for the proofs.
Since $\bm M_c$ is invertible, the acceleration of the dependent generalized coordinates can be always computed from  \eqref{eq_EqMotion2} regardless the ill-conditioning of the constraint, i.e.,
\begin{align} \notag
\ddot{\bm q} & = -\bm M_c^{-1} \bm P (\bm u -\bm h) + \big( \nu \bm M_c^{-1} \bm P\bm\Omega -  \bm M_c^{-1}\bm P \bm M \bm\Omega  \big)  \dot{\bm q}  \\ \label{eq:ddq_sim}
& = \bm M_c^{-1} \bm P (\bm u - \bm h) + \bm S \bm\Omega \dot{\bm q}
\end{align}
Now suppose vector $\ddot{\bm q}^* = \bm M^{-1} (\bm u -\bm h)$ virtually represent the acceleration of the MBS when all constraints are ignored. Then, by virtue of definition of the unconstrained acceleration and Properties~\ref{prop:Mc}, one can equivalently rewrite \eqref{eq:ddq_sim} in the following form
\begin{equation} \label{eq:ddq_S}
\ddot{\bm q} = (\bm I - \bm S) \ddot{\bm q}^*  + \bm S \bm\Omega \dot{\bm q}.
\end{equation}
\begin{remark}
Equation \eqref{eq:ddq_S} reveals that the acceleration vectors of constrained and unconstrained MBSs are related by the oblique projection matrix $\bm I- \bm S$.
\end{remark}
Notice that the generalized velocities are not continuous at impacts and therefore the post-impact velocities can not determined from  \eqref{eq:ddq_S}. Nevertheless, it will be shown in the next section \eqref{sec:impact} that the post-impact velocity can be also computed from the oblique projection.

\subsection{Constraint force and oblique projection}
Upon substitution of the acceleration from  \eqref{eq:ddq_sim} into \eqref{eq:Mddotq_standard} and rearranging the latter equation, we arrive at the equation the constraint forces in the following compact form
\begin{equation} \label{eq:fc_R}
\bm f =  \bm S^T \big( \bm h - \bm u \big) + \bm S^T \bm M \bm\Omega \dot{\bm q}
\end{equation}
Since the projection matrix $\bm S$ is always well-defined, equation \eqref{eq:fc_R} gives a unique solution to the constraint forces. However, the generalized Lagrangian multipliers cannot be determined uniquely unless the constraints are linearly independent. In any case, the pseudo-inverse can be used for a minimum-norm solution, i.e., $\bm\lambda = \bm A^{+T} \bm f$.

\section{Passive unilateral constraint becomes active: Generalized Impact} \label{sec:impact}
When a passive unilateral constraint becomes active such event is usually accompanied by corresponding negative velocity $\dot{\phi}_{u_i} < 0$. This event  causes impact accompanied by impulsive forces and sudden change in the generalized velocity and therefore a discontinuity in the velocities. Due to discontinuity of the velocities during impact event,  the velocity  instantaneously jumps from  $\dot {\bm q}^-$ to $\dot{\bm q}^+$ and thus requiring infinitely large acceleration and constraint force. Therefore, the acceleration model \eqref{eq:ddq_S} is not adequate to determine post-impact velocities rather an impact model is required to deal with the impulsive constraint force and discontinuities in the velocities. This section presents a  closed-form solution to the generalized impulse-momentum equation of constrained MBSs by establishing explicit relationship between the pre- and post-impact velocities based on the oblique projection $\bm S$.

Although it is not possible to calculate the velocity change at the time of impact through integration of the equations of motion, it is possible to calculate the velocity change using the {\em Newton's impact law}. Suppose multiple impacts  occur at time interval $[t^-, \; t^+]$ and the impact duration $\delta t= t^+ - t^-$ is infinitesimal. Also, let us define the $m_u  \times m_u$ switching matrix $\bm\Gamma^*$ whose corresponding diagonal entries switch from zero to one upon detection of impact, i.e.,
\begin{equation} \label{eq:Gamma*}
\Gamma^*_{ii} =\left\{ \begin{array}{ll} 1 & \mbox{if} \qquad \phi_i \leq 0 \quad \wedge \quad \dot{\phi}_i < 0
\\ 0 & \mbox{otherwise} \end{array} \right.
\end{equation}
Notice that simultaneous activation of  more than one unilateral constraint is possible and therefore matrix $\Gamma _{ii}$ may contain more than one nonzero entries. For instance, if occurrence of simultaneous activation of the first and third  unilateral constraints are detected, then the switching matrix becomes
\begin{equation}
\bm\Gamma^*  = \mbox{diag}\big\{ [1, 0, 1, 0, \cdots, 0] \big\}
\end{equation}
From the impact switching matrix, it follows the jacobian matrix during the impact as
\begin{equation} \label{eq:Adef_star}
\bm A_u  = (\bm\Gamma + \bm\Gamma^* ) \frac{\partial \bm\phi_u}{\partial \bm q},
\end{equation}
and subsequently \eqref{eq:Adef}  can be treated as the mapping from the lagrangian multipliers associated with both bilateral constraints and the activated unilateral constraints to the joint torques. In the followings, we assume the overall jacobian matrix $\bm A$ accounts for the simultaneously activated  unilateral constraints according to \eqref{eq:Adef_star}. Since the generalized coordinate $\bm q$ is constant during the impact, the mass matrix $\bm M(\bm q)$ and the Jacobian $\bm A (\bm q)$ remain unchanged during  the impact. Therefore, one can carry out integration of the differential equation \eqref{eq:Mddotq_standard} over $[t^-, \; t_+]$ to obtain impact equation of the system as
\begin{equation} \label{eq:Dirac}
\bm M (\dot{\bm q}^+ - \dot{\bm q}^- ) = \bm i_u + \bm A^{T}  \bm i_\lambda.
\end{equation}
Here, $\bm i_u$ and $\bm i_\lambda$ are the impacts or the Dirac integrals of the input forces and lagrangian multipliers, i.e.,
\begin{equation}
\bm i_u  = \lim_{\delta t \rightarrow 0} \int_{t^-}^{t^- + \delta t} \bm u \; {\rm d}t, \qquad  \bm i_\lambda  = \lim_{\delta t \rightarrow 0} \int_{t^-}^{t^- + \delta t} \bm\lambda \; {\rm d}t,
\end{equation}
and $\dot{\bm q}^- = \dot{\bm q}(t^-)$ and $\dot{\bm q}^+ = \dot{\bm q}(t^+)$ are pre-impact and post-impact velocities.  Notice that in derivation of \eqref{eq:Dirac}, we also assumed  $\bm h$ to be a continuous function and therefore its integration over  $[t^-, \; t^+]$ vanishes. Now suppose $\bm P $ be the projection matrix associated with the impact Jacobian matrix. Then, by pre-multiplying both sides of \eqref{eq:Dirac} by $\bm P $, we arrive at
\begin{equation} \label{eq:P*M}
\bm P  \bm M ( \dot{\bm q}^+ -\dot{\bm q}^-)= \bm P  \bm i_u
\end{equation}
The above equation  does not yet completely describe the impact, because the impact
laws are still missing. The above equation constitutes the balance of momentum in the contact that in conjunction with a  restitution law can uniquely determine the relationship between the pre-impact and post-impact velocities. There are different approaches to define the restitution coefficient based on either kinematic treatment or energy considerations of the impact \cite{Khulief-2013}.

There are three main definitions for the restitution coefficient. They are: i) Newton's definition  (Kinematic), ii) Poisson's definition (Kinetic), and iii) Stronge's definitions (Energetic) \cite{Ahmad-Ismail-Mat-2016}. The pros and  cons of restitution coefficient models based on these definitions have been evaluated and a survey can found, for example, in \cite{Ahmad-Ismail-Mat-2016}. It is known that all these definitions of coefficient of restitution are equivalent for frictionless contact, albeit the Stronge's definition might result in a better solution in terms of energy-inconsistency of the post-impact if the configuration of impact involves with friction and the direction of the slip changes \cite{Ismail-Stronge-2008}. Nevertheless, the {\em Newton}'s rule is the most widely used restitution rule and it is reasonably accurate for bulky bodies \cite{Khulief-2013}. The Newton's restitution law is a kinematic form which characterizes the normal component of the velocities of two rigid bodies before and after impact by a kinematic parameter called the restitution coefficient, i.e.,
\begin{equation} \label{eq:restitution_law}
\bm\Gamma^*  \dot{\bm \phi}_{u}^+ = -e \bm\Gamma^*  \dot{\bm\phi}_{u}^-,
\end{equation}
where $0<e<1$ is the global {\em restitution coefficient}. The assumption that all the multiple contacts have identical restitution coefficients will be later relaxed in Section~\ref{sec:local_restitution_coeff}. The restitution equations \eqref{eq:restitution_law} can be equivalently written in terms of the generalized velocities before and after the impact by
\begin{equation} \label{eq:GammaAdq}
\bm\Gamma^*  \bm A_u  \dot{\bm q}^+ = -e \bm\Gamma^*  \bm A_u \dot{\bm q}^-
\end{equation}
On the other hand, admissible pre- and post-impact velocities should remain in the null-space of the Jacobian matrix \eqref{eq:Adef}, i.e., $\bm\Gamma  \bm A_u  \dot{\bm q}^+ = \bm\Gamma  \bm A_u  \dot{\bm q}^- = \bm 0$. Thus, \eqref{eq:GammaAdq} can be equivalently written as
\begin{equation} \label{eq:A*dq}
\bm A  \dot{\bm q}^+ = - e  \bm A  \dot{\bm q}^- \quad \mbox{or} \quad \bm A  \big( \dot{\bm q}^+ + e \dot{\bm q}^- \big) = \bm 0,
\end{equation}
which implies vector $\dot{\bm q}^+ + e \dot{\bm q}^-$ should belong to the null-space of the impact jacobian matrix $\bm A $. Using the notion of the projection matrix, \eqref{eq:A*dq} can be equivalently written as follow
\begin{equation} \label{eq:I-P*}
( \bm I - \bm P ) \dot{\bm q}^+ = -e (\bm I- \bm P  )\dot{\bm q}^-.
\end{equation}
Substituting  expression $\dot{\bm q}^+ = \bm P  \dot{\bm q}^+ - e \dot{\bm q}^- + e \bm P  \dot{\bm q}^-$  obtained from  \eqref{eq:I-P*} into \eqref{eq:P*M} gives
\begin{equation} \label{eq:Mo_dotq+}
\bm P  \bm M \bm P  \dot{\bm q}^+ = (e + 1) \bm P  \bm M \dot{\bm q}^- - e \bm P  \bm M \bm P  \dot{\bm q}^- + \bm P  \bm i_u
\end{equation}
Equations \eqref{eq:Mo_dotq+} together with  \eqref{eq:I-P*} complete the impact problem because they provide us with sufficient equations to uniquely determine the post-impact velocity from the pre-impact velocity. Pre-multiplying \eqref{eq:I-P*} by virtual mass $\nu$ and adding both sides of the resulting equation with those of  \eqref{eq:Mo_dotq+}, we obtain
\begin{equation}
{\bm M} _c \dot{\bm q}^+ = (e +1) \bm P  \bm M \dot{\bm q}^- - e {\bm M} _c \dot{\bm q}^- +  \bm M_c^{-1} \bm P  \bm i_u
\end{equation}
where ${\bm M} _c$ is the  constraint mass matrix.  We are now able to propose the projection-based generalized impact model of multibody system described in  \eqref{eq:Dirac} and \eqref{eq:restitution_law}. Since ${\bm M} _c$ is always invertible, the above equation can be solved  through matrix inversion using identity \eqref{eq:(I-S)invM} as
\begin{align} \notag
\dot{\bm q}^+ &= \dot{\bm q}^- -(e +1) \bm S  \dot{\bm q}^- +  \bm M^{-1}(\bm I - \bm S^T ) \bm i_u\\\label{eq:dot_q_plus}
& = \dot{\bm q}^- -(e +1) \bm S  \dot{\bm q}^- + \bm M_o^+ \bm i_u,
\end{align}

\begin{remark}
It can be inferred from \eqref{eq:dot_q_plus} that the external or actuation impulse has no effect on the post-impact velocity if $\bm i_u \in \mathcal{R}(\bm S^T)$.
\end{remark}

It worths mentioning that  $\bm i_u$  represents all external impulse acting on the MBSs including those from the actuators. Most actuators can not respond quickly enough to generate impacts, albeit{\em impact actuators} also do exist.

To this end, it should be also mentioned that  modelling the coefficient of restitution  has been the subject of intensive investigation by many researchers using theoretical, numerical, and experimental methods \cite{Kangur-Kleis-1988,Zhang-Vu-Quoc-2002,Najafabadi-Kovecses-Angeles-2008,Jackson-Green-2009}. It is well known that the coefficient of restitution changes with not only  material properties such as yield strength, elastic modulus, density, and Poisson's ratio, but also with pre-impact velocity \cite{Lankarani-Nikravesh-1990,Marhefka-Orin-1999,Zhang-Vu-Quoc-2002,Najafabadi-Kovecses-Angeles-2008,Jackson-Green-2009}. However, this is not a limitation of the momentum-balance and restitution law method because the coefficient of restitution  can be predicted as a function of the pre-impact velocity when the specifics of the collision are simplified.

\subsection{Geometrical interpretation of generalized impact in the joint space}
It can be also readily inferred from \eqref{eq:dot_q_plus} that for the case of perfectly elastic impact, i.e., $e =1$, and no external impulse, i.e., $\bm i_u = \bm 0$,  the relation \eqref{eq:dot_q_plus} simply becomes
\begin{subequations} \label{eq:post-dq}
\begin{equation} \label{eq:Rdq+}
\dot{\bm q}^+ = \bm R   \dot{\bm q}^- \qquad \quad \leftarrow \quad e=1
\end{equation}
where $\bm R =\bm I - 2 \bm S$  and
\begin{equation}
\bm R^{2} = \bm I.
\end{equation}
Therefore, the transformation matrix $\bm R $ is {\em involuntary} meaning that the transformation is its own inverse.  In other words, $\bm R $ is a {\em reflection} matrix written in terms of projections. This is interesting results because \eqref{eq:Rdq+} represents n-dimensional extension of elastic impact of a simple point mass scenario. On the other hand, in the case of contact with zero-restitution $e=0$, the pos-impact velocity and pre-impact velocity are related by the oblique projection operator  $\bm I - \bm S $, i.e.,
\begin{equation} \label{eq:Rdq+2}
\dot{\bm q}^+ = (\bm I - \bm S )  \dot{\bm q}^- \qquad \quad \leftarrow \quad e=0
\end{equation}
In practice, when $e\in[0 \; 1]$, the pos-impact velocity can be obtained by linear  combination of equations \eqref{eq:Rdq+} and \eqref{eq:Rdq+2} as follow
\begin{equation} \label{eq:impact_geometry}
\dot{\bm q}^+ = \Big( e \bm R  +(1 - e)(\bm I - \bm S ) \Big) \dot{\bm q}^-.
\end{equation}
\end{subequations}
The geometrical interpretation of the projection-based impact model is given as follow
\begin{remark}
The pre-impact and post-impact velocities of a MBS with chain, tree and close-loop topologies
are related by the reflection matrix $\bm R $ when $e=1$ and by the oblique projection matrix $\bm I - \bm S $ when $e=0$,  and by a linear combination of the two transformation matrices when $e \in [0, \; 1]$.
\end{remark}

Defining variables $\bm p^- = \bm M \dot{\bm q}^-$ and $\bm p^+ = \bm M \dot{\bm q}^+$ as the  pre-impact and post-impact generalized momentum of the MBS and using identity \eqref{MS=StM} in \eqref{eq:post-dq}, the pre- and post-impact states of a MBS can be alternatively described in terms of the generalized momentums as follow
\begin{align} \notag
{\bm p}^+ &= {\bm p}^- -(e +1)\bm S^{T} {\bm p}^- \qquad \forall e \in[0, \; 1] \\ \label{eq:momentum_geometry}
& =\Big( e \bm R^{T} +(1 - e)(\bm I - \bm S^{T}) \Big) {\bm p}^-
\end{align}

\subsection{Impact calculation}
Suppose  $\bm i_f= \bm A^{T} \bm i_{\lambda}$ denotes the generalized impact in MBS. Then substituting the velocity change $\dot{\bm q}^+ -  \dot{\bm q}^-$ from \eqref{eq:dot_q_plus} into \eqref{eq:Dirac} and using \eqref{MS=StM}, we arrive at the expression of the generalized impact
\begin{align} \notag
\bm i_f  & = -(e + 1) \bm M \bm S \dot{\bm q}^- - \bm S^{T} \bm i_u\\ \label{eq:bar_fc}
& = -(e + 1) \bm S^{T} \bm p^- - \bm S^{T} \bm i_u
\end{align}
\begin{remark}
The input impact, $\bm i_u$, has not effect on the generalized impact, $\bm i_f$, if $\bm i_u \in {\cal R}(\bm I - \bm S^T)$. Moreover,  the transition from non-contact to contact occurs with zero impact  if ${\bm p}^-, \bm i_u \in {\cal R}(\bm I - \bm S^T)$.
\end{remark}
\begin{remark}
The generalized impact force in a MBS with chain, tree and close-loop topologies is proportional to the  projected version of its pre-impact generalized momentum.
\end{remark}
Notice that $\bm i_f$ can be always uniquely determined from \eqref{eq:bar_fc} but that is not the case for $\bm i_{\lambda}$.  One can conclude from \eqref{(I-P)St} and \eqref{eq:SP=0} that ${\cal N}^{\perp}(\bm A) \equiv {\cal R}(\bm S)$ or equivalently ${\cal R}(\bm S^T) \equiv {\cal R}(\bm A^T)$. Therefore,  the solution obtained from \eqref{eq:bar_fc} must be in  ${\cal R}(\bm A^{T})$, which  means that at least  one solution to equation $\bm i_f = \bm A^T \bm i_{\lambda} $ is guaranteed and that can be obtained through pseudo-inversion $\bm i_{\lambda} = \bm A^{+T} \bm i_f$.

\subsection{Energetic consistency}
Energy lost during an impact is an important quantity not only to gain insight into complex physical phenomenon during contact but to examine whether an impact model is physically consistent \cite{Stronge-1991,Stoianovici-Hurmuzlu-1996,Najafabadi-Kovecses-Angeles-2008}. Suppose $K^-=\frac{1}{2} \dot{\bm q}^{-T}  \bm M \dot{\bm q}^-$ and $K^+=\frac{1}{2} \dot{\bm q}^{+T}  \bm M \dot{\bm q}^+$ are the pre-impact and post-impact values of the kinetic energy of the constrained system. Then, by virtue of identities \eqref{eq:dot_q_plus} and \eqref{MS=StM} one can derive the expression of the difference between the post-impact  and pre-impact kinetic energy as
\begin{equation} \label{eq:delta_K}
 K^+ - K^- = -\frac{1}{2}(1-e^2) \dot{\bm q}^{-T} \bm S^{T} \bm M \bm S  \dot{\bm q}^- + \frac{1}{2} \bm i_u^T \bm M_o^+  \bm i_u + \dot{\bm q}^{-T} (\bm I - \bm S^T) \bm i_u,
\end{equation}
which is in the form of a quadratic function of the pre-impact velocity, $\dot{\bm q}^-$, and the external impact, $\bm i_u$, variables. In the absence of external impact, the energy loss absorbed in the contact is $W_{\rm loss}=K^+ - K^-$ with $\bm i_u \equiv \bm 0$, i.e.,
\begin{equation} \label{eq:kinetic}
W_{\rm loss} = - \frac{1}{2} (1-e^2) \dot{\bm q}^{-T} \bm S^{T} \bm M  \bm S \dot{\bm q}^-
\end{equation}
where $\bm S^{T} \bm M  \bm S$ is called {\em projected inertia matrix}.
\begin{properties}
The projected inertia matrix satisfies
\begin{subequations}
\begin{align} \label{eq:SMT_psd}
& \bm S^{T} \bm M  \bm S \succeq  0 \\ \label{eq:M-SMT_psd}
& \bm M - \bm S^{T} \bm M  \bm S  \succeq 0
\end{align}
\end{subequations}
\end{properties}
Inequality \eqref{eq:SMT_psd} can be readily inferred from that $\bm M \succeq  0$. Moreover by virtue of \eqref{eq:PMP}, we have
\begin{align*}
\bm M -  \bm S^{T} \bm M  \bm S & =  \bm M -   \bm M  \bm S \\
& = \bm M - \bm M + \bm M \bm P \bm M_c^{-1} \bm M\\
&= \bm M \bm P \bm M_c^{-1} \bm P \bm M \preceq 0
\end{align*}
which is a positive semi-definite matrix. In other words, \eqref{eq:M-SMT_psd} implies
\begin{equation} \label{eq:SMS<M}
\bm\zeta^T \bm S^{T} \bm M \bm S \bm\zeta  \leq \bm\zeta^T \bm M \bm\zeta  \qquad \forall \bm\zeta \in \mathbb{R}^n,
\end{equation}
which proves \eqref{eq:M-SMT_psd}. Equation \eqref{eq:kinetic} reveals that the energy absorbed in the impact is a quadratic function of the pre-impact velocity is proportional to $1-e^2$, hence the energy loss becomes zero when $e=1$. Now, let us define quantity
\begin{align} \notag
\gamma &= \frac{K^+}{K^-} =  1 + \frac{W_{\rm loss}}{K^-} \\ \label{eq:beta}
&= 1 - (1- e^2) \frac{\dot{\bm q}^{-T} \bm S^{T} \bm M  \bm S \dot{\bm q}^-}{\dot{\bm q}^{-T}\bm M  \dot{\bm q}^- }
\end{align}
which measures energy lost during an impact. For instance, in the case of perfect elastic impact where $e=1$ when $\gamma=1$, i.e., no energy loss. From  \eqref{eq:SMT_psd} and \eqref{eq:SMS<M}, one can infer that $\gamma$ should remain in the following bound
\begin{equation}
0 \leq \gamma \leq 1  \quad \leftarrow \quad e \in[0, \; 1].
\end{equation}
Thus
\begin{equation} \notag
W_{\rm loss} \leq 0 \quad \mbox{and} \quad | W_{\rm loss} | \leq K^-
\end{equation}
In other words, the post-impact energy of the mechanical system is less than or equal to the pre-impact energy level meaning that the projection-based impact law \eqref{eq:impact_geometry} satisfies the fundamental energetic consistency.

\begin{remark}
It is evident from \eqref{eq:SMS<M} that the impact model for a MBS with chain, tree and close-loop topologies is physically consistent, i.e., post-impact kinetic energy never exceeds  that of pre-impact, if the global restitution coefficients $e$ is bounded within $[0, 1]$ and that the energy at impacts is conserved if $e=1$.
\end{remark}

It is worth noting from expression \eqref{eq:beta}  that $\gamma$ is affected not only
by the coefficient of restitution, but also the configuration of
the whole multibody system and the direction of the velocity vector. The energy dissipation done by the contact force $\bm f$ along the complete path of compression and restitution can be expressed as
\begin{align} \notag
W_{\rm loss} &= \oint \bm f \cdot d \bm q = \int_t^{t + \delta t} (\bm f \cdot \dot{\bm q}) dt \\
&  = - \frac{1}{2} (1 - e^2) \dot{\bm q}^{-T} \bm S^{T} \bm M \bm S  \dot{\bm q}^-
\end{align}
in which the last term is concluded from \eqref{eq:kinetic}. On the other hand, from expressions \eqref{eq:dot_q_plus} and \eqref{eq:bar_fc}, one can infer
\begin{equation}
{\bm i}_f \cdot  \bar{\dot {\bm q}} = - \frac{1}{2} (1 - e^2) \dot{\bm q}^{-T} \bm S^{T} \bm M \bm S  \dot{\bm q}^- \quad \mbox{where} \quad  \bar{\dot {\bm q}} = \frac{1}{2}(\dot{\bm q}^- + \dot{\bm q}^+)
\end{equation}
is defined as the average of pre- and post-impact velocities. Thus
\begin{equation}
W_{\rm loss} = \int_t^{t + \delta t} (\bm f \cdot \dot{\bm q}) dt = {\bm i}_f \cdot  \bar{\dot {\bm q}}
\end{equation}
Clearly, if the impact is lossless, then the inner product must be zero meaning that ${\bm i}_f \perp  \bar{\dot {\bm q}}$. Defining $\bar{\dot {\bm\phi}}_ = \bm A  \bar{\dot{\bm q}}$, one can verify that the following identity
\begin{equation} \label{eq:ifdotphi}
\bm i_f \cdot \bar{\dot{\bm q}} = \bm i_{\lambda} \cdot \bar{\dot {\bm\phi}} = \bm i_{\lambda_u} \cdot  \bar{\dot {\bm\phi}}_u
\end{equation}
which indicate the work done by the contact at the joint space and the task space are equal. Note that the latter equality in \eqref{eq:ifdotphi} is obtained from the fact that the bilateral constraints do not perform any work because $\dot {\bm\phi}_b \equiv \bm 0$. Moreover from the restitution law we have $\bar{\dot{\bm\phi}}_u = \frac{1}{2}(1 - e)\dot{\bm\phi}_u^-$ and therefore we can say
\begin{equation} \label{eq:Wloss/1-e}
\bm i_{\lambda_u} \cdot  {\dot {\bm\phi}}_u^- = \frac{2 W_{\rm loss}}{1 - e} =-(1+e) \dot{\bm q}^{-T} \bm S^{T} \bm M \bm S  \dot{\bm q}^- \leq 0
\end{equation}

\begin{remark}
For the case of $e=1$, the kinetic energy ratio is one meaning that the kinetic energy is preserved regardless the direction of the pre-impact velocity vector or configuration of the whole multibody system.
\end{remark}

\section{Local coefficients of restitution} \label{sec:local_restitution_coeff}
The projection-based formulation of the generalized impact law presented in the previous section assumes a global restitution coefficient for all contacts. However, it is known that the restitution coefficient depends on many factors such as the materials and
the pre-impact velocity  \cite{Najafabadi-Kovecses-Angeles-2008}. More specifically, the dependency of the coefficient of restitution to the surface topography in addition to the material properties, e.g., Young's modulus-to yield stress ratio, and incident velocity is shown in \cite{Lu-Kuo-2003,Khulief-2013}.
Therefore, one may not be able to assume  global coefficient of restitution to develop a multiple impact model, rather local coefficients of restitution have to be incorporated in the model. This section prolongs the projection-based model presented in the previous section for multiple simulations contacts with non-identical coefficients of restitution.

The restitution rule \eqref{eq:restitution_law}-\eqref{eq:Wloss/1-e} for the case of non-identical coefficients of restitution can be extended to the following matrix form
\begin{equation} \label{eq:dq+matrix}
\bm A  \dot{\bm q}^+ = - \bm E \bm A  \dot{\bm q}^- ,
\end{equation}
where matrix $\bm E$ contains local coefficients of restitution. For instance, the restitution matrix may take the form
\begin{equation}
\bm E = \mbox{diag}\big(0, \cdots, 0, e_1, e_1, \cdots, e_m \big),
\end{equation}
where $e_i$ represents the restitution coefficient at the $i$th contact as the zeros are associated with the bilateral constraints. In order to be able to represent the impact formulation in a closed-form, it becomes necessary to transform the restitution matrix into the following form
\begin{equation}
\tilde{\bm E} = \bm A^{+} \bm E \bm A
\end{equation}
Notice that even if matrix $\bm E$ is chosen to be diagonal, matrix $\tilde{\bm E}$ is not necessarily diagonal. From the above definition, one can really verify the following identities
\begin{equation} \label{eq:alphaP}
\tilde{\bm E}  \bm P  =  \bm P  \tilde{\bm E} = \bm 0
\end{equation}
Pre-multiplying both sides of \eqref{eq:dq+matrix} by $\bm A  \bm A^{+}$ yields
\begin{equation} \label{eq:A2dq=0}
\bm A  \dot{\bm q}^+ =  -\bm A  \tilde{\bm E} \dot{\bm q}^- , \quad \mbox{or} \quad  \bm A  (\dot{\bm q}^+ + \tilde{\bm E} \dot{\bm q}^-) = \bm 0,
\end{equation}
in which we use the property of pseudo-inverse $\bm A \bm A^+ \bm A = \bm A$. In view of \eqref{eq:alphaP}, the above equation can be equivalently written in terms of the corresponding projection matrix as follow
\begin{align}
( \bm I - \bm P ) \dot{\bm q}^+ & = -\tilde{\bm E} (\bm I- \bm P  )\dot{\bm q}^- \\ \label{eq:I-P*2}
 & = -(\bm I- \bm P  ) \tilde{\bm E}   \dot{\bm q}^-,
\end{align}
which resembles \eqref{eq:I-P*}. Finally, in a development similar to \eqref{eq:I-P*}-\eqref{eq:dot_q_plus}, one can combine the kinematic equation \eqref{eq:I-P*2} with the momentum balance equation \eqref{eq:P*M} to
arrive at the following expression
\begin{equation} \notag
\dot{\bm q}^+ = \dot{\bm q}^- - \bm S  (\tilde{\bm E} + \bm I) \dot{\bm q}^- +\bm M_o^+ \bm i_u.
\end{equation}
If $\bm i_u \equiv 0$, then the above equation becomes
\begin{equation} \label{eq:dot_q_plus2}
\dot{\bm q}^+ = \big( \bm R  \tilde{\bm E}  +( \bm I - \bm S ) (\bm I - \tilde{\bm E}) \big) \dot{\bm q}^-.
\end{equation}
Subsequently, the corresponding impulse can be found by substituting $\dot{\bm q}^+$ from \eqref{eq:dot_q_plus2} into \eqref{eq:Dirac}
\begin{equation} \label{eq:if2}
\bm i_f = - \bm M \bm S (\tilde{\bm E} + \bm I) \dot{\bm q}^-
\end{equation}
Equations \eqref{eq:dot_q_plus2} and \eqref{eq:if2} constitute  the generalized impact model of multibody system involving simultaneous multiple impacts with non-identical coefficients of restitution.

Similar to \eqref{eq:kinetic}, the energy loss can be algebraically obtained  from the expression of the post-impact velocity \eqref{eq:dot_q_plus2} as follow
\begin{align} \notag
 & W_{\rm loss} = -\frac{1}{2}\dot{\bm q}^{-T}  \tilde{\bm E}^T \bm S^{T} \bm M \dot{\bm q}^- -   \frac{1}{2}\dot{\bm q}^{-T}   \bm S^{T} \bm M \dot{\bm q}^- \\ \notag
 & - \frac{1}{2}\dot{\bm q}^{-T}   \bm M \bm S\dot{\bm q}^- - \frac{1}{2}\dot{\bm q}^{-T}   \bm M \dot{\bm q}^- + \frac{1}{2}\dot{\bm q}^{-T}  \tilde{\bm E}^T \bm S^{T} \bm M \bm S \tilde{\bm E} \dot{\bm q}^- \\ \notag
 & + \frac{1}{2}\dot{\bm q}^{-T}  \tilde{\bm E}^T \bm S^{T} \bm M \bm S  \dot{\bm q}^- +  \frac{1}{2}\dot{\bm q}^{-T}  \tilde{\bm E}^T \bm S^{T} \bm M \bm S \tilde{\bm E} \dot{\bm q}^-  + \frac{1}{2}\dot{\bm q}^{-T}   \bm M \dot{\bm q}^- \\ \label{eq:Wloss_E1}
& =  \frac{1}{2}\dot{\bm q}^{-T}  \tilde{\bm E}^T \bm S^{T} \bm M \bm S \tilde{\bm E} \dot{\bm q}^- - \frac{1}{2}\dot{\bm q}^{-T}   \bm S^{T} \bm M \bm S \dot{\bm q}^-  \\ \label{eq:Wloss_E2}
& = - \frac{1}{2}\dot{\bm q}^{-T} (\bm I - \tilde{\bm E}^T) \bm S^{T} \bm M  \bm S (\bm I + \tilde{\bm E})\dot{\bm q}^-
\end{align}
Using expression  \eqref{eq:dot_q_plus2} in a development similar to \eqref{eq:kinetic}-\eqref{eq:beta}, one can derive the post-impact to pre-impact energy ratio  for the case of non-identical coefficients of restitution as
\begin{equation}  \notag
\gamma = 1 - \frac{\bm\zeta^{T} (\bm I - \tilde{\bm E}^T) \bm S^{T} \bm M  \bm S (\bm I + \tilde{\bm E})\bm\zeta} {\bm\zeta^T \bm M   \bm\zeta }.
\end{equation}
The expression of energy loss in \eqref{eq:Wloss_E1} can be also equivalently written in the following compact form by changing the variable from $\dot{\bm q}^-$ to $\dot{\bm\phi}^-=\bm A  \dot{\bm q}^-$
\begin{equation} \label{eq:loss2}
W_{\rm loss} = \frac{1}{2} \dot{\bm\phi}^{-T} \big( \bm E\bm Q \bm E - \bm Q \big) \dot{\bm\phi}^{-},
\end{equation}
where $\bm Q=\bm Q^T \succeq 0$ is a semi-positive matrix defined as follow
\begin{equation} \label{eq:QMI_1}
\bm Q = \bm G^T \bm M \bm G,  \quad \mbox{and} \quad  \bm G =\bm S  \bm A^{+}.
\end{equation}
Clearly, energy loss during the impact is non-positive if the matrix $\bm E\bm Q \bm E - \bm Q$ is semi-negative definite. However, the latter is not generally obvious unless for spacial cases. For instance, if $\bm E=e \bm I$, the above matrix expression comes down to $(e^2-1)\bm Q$ and if $e\leq 1$ then the matrix  is automatically semi-negative definite because $\bm Q \succeq 0$. In general, we can say the multiple impact model for the case of non-identical restitution coefficients is energetically consistent if the restitution coefficient matrix $\bm E$ satisfies the following {\em quadratic matrix inequality} (QMI)
\begin{equation} \label{eq:QMI}
\bm E\bm Q \bm E - \bm Q \preceq 0  \qquad \mbox{where} \quad \bm E, \bm Q \succeq 0,
\end{equation}
in which matrix $\bm E$ is the variable. The QMI \eqref{eq:QMI} in $\bm E$ is simplified version of generalized Riccati inequality and therefore it can be also expressed as the following {\em linear matrix inequality} (LMI) in $\bm E$ by applying the Schur complement lemma \cite{VanAntwerp-Braatz-2000}:
\begin{equation} \label{eq:LMI}
\begin{bmatrix} \bm G^T \bm M \bm G & \bm E \bm G^T \\ \bm G \bm E & \bm M^{-1} \end{bmatrix} \succeq 0
\end{equation}
A number of toolboxes are available \cite{Gahinet-Nemirovski-Laub-1995,Delebecque-Nikoukhah-1995} for solving the above LMI problem to find the feasible set of restitution matrix $\bm E$ that results that makes the system energetically consistent.

\appendix

\subsection{Time-derivative of projection matrix} \label{appx_skew}
The Tikhonov regularization theorem \cite{Golub-VanLoan-1996} describes the pseudo-inverse as  the following limit
\begin{equation} \label{eq:Tikhonov}
\bm A^+ = \lim_{\epsilon \rightarrow 0} \bm A^T(\bm A \bm A^T + \epsilon \bm I)^{-1}
\end{equation}
By differentiation of  the above expression, one can verify that   the time-derivative of the  pseudo-inverse can written in the following form
\begin{align} \notag
\frac{d}{dt} \bm A^+ & =\lim_{\epsilon \rightarrow 0} \dot{\bm A}^T(\bm A \bm A^T + \epsilon \bm I)^{-1} - \bm A^T (\bm A \bm A^T + \epsilon \bm I)^{-1}  [\dot{\bm A} \bm A^T + \bm A \dot{\bm A}^T]  (\bm A \bm A^T + \epsilon \bm I)^{-1} \\ \notag
& = - \bm A^+ \dot{\bm A} \bm A^+ + \lim_{\epsilon \rightarrow 0} \dot{\bm A}^T(\bm A \bm A^T + \epsilon \bm I)^{-1}  - \bm A^+ \bm A \dot{\bm A}^T(\bm A \bm A^T + \epsilon \bm I)^{-1}\\\label{eq:dotA+}
&= -\bm A^+ \dot{\bm A} \bm A^{+} + \lim_{\epsilon \rightarrow 0} \bm P \dot{\bm A}^T (\bm A \bm A^T + \epsilon \bm I)^{-1}
\end{align}
On the other hand, using \eqref{eq:dotA+} in the time-derivative of the expression of the projection matrix $\bm P= \bm I - \bm A^+ \bm A$ yields
\begin{equation} \label{eq:dotP-A+}
\begin{split}
\dot{\bm P} &= - \frac{d}{dt}\bm A^+ \bm A - \bm A^+ \dot{\bm A} \\ \notag
& =  \bm A^+ \dot{\bm A} \bm A^{+} \bm A + \lim_{\epsilon \rightarrow 0} \bm P \dot{\bm A}^T (\bm A \bm A^T + \epsilon \bm I)^{-1} \bm A - \bm A^+ \dot{\bm A} \\
& =\bm A^+ \dot{\bm A} (\bm I - \bm P) + \bm P \dot{\bm A}^T \bm A^{+T} - \bm A^+ \dot{\bm A} \\
&=  \bm\Lambda + \bm\Lambda^T
\end{split}
\end{equation}
where $\bm\Lambda=- \bm A^+ \dot {\bm A} \bm P$. Note that identity $\bm P \bm A^+ = \bm A^{+T} \bm P = \bm 0$  implies that $\bm\Lambda^T \bm P = \bm 0$ and hence one can conclude $\ddot{\bm q}_{\perp}= \dot{\bm P} \dot{\bm q}= \bm\Lambda \dot{\bm q} = \bm\Omega \dot{\bm q}$ \cite{Aghili-2015b}.

\subsection{Properties of $\bm M_c$} \label{appx:Mc_properties}
Consider non-zero vector
$\bm a \in \mathbb{R}^n$ and its orthogonal decomposition components $\bm a_{\parallel}=\bm P \bm a$ and $\bm a_{\perp}=(\bm I - \bm P) \bm a$. Then, one can say
\begin{equation} \label{eq:twoMass}
\bm a^T {\bm M}_c \bm a  =
\bm a^T_{\parallel} \bm M \bm a_{\parallel} + \nu \| \bm a_{\perp}
\|^2
>0,
\end{equation}
Notice that both terms $\bm a^T_{\parallel} \bm M \bm a_{\parallel}>0$ and $\| \bm a_{\perp}
\|^2>0$ are positive semi-definite. Moreover, for a given non-zero vector $\bm a$ if $\bm a_{\perp}= \bm 0$ then $\bm a^T_{\parallel} \neq 0$ and vice versa. This means that the summation of the two orthogonal terms must be positive definite and so must be the constraint inertia matrix ${\bm M}_c$.

By definition we have
\begin{equation} \notag
\bm\Omega \dot{\bm q} = (\bm I - \bm P) \bm\Omega \dot{\bm q}
\end{equation}
Since $\bm M_c$ is always invertible, i.e., $\bm M_c^{-1}\bm M_c = \bm I$,  the above equation can be equivalently written as
\begin{align*}
\bm\Omega \dot{\bm q} & = \bm M_c^{-1} \bm M_c (\bm I - \bm P) \bm\Omega \dot{\bm q} \\
& = \bm M_c^{-1} \nu (\bm I - \bm P)  \bm\Omega \dot{\bm q} \\
& = \nu \bm M_c^{-1}\bm\Omega \dot{\bm q}
\end{align*}
which proofs \eqref{eq:McinvOmega}.

From definition we have $\bm M_c \bm P = \bm M$ and hence $\bm M \bm M_c^{-1} = \bm M_c \bm P \bm M_c^{-1} = \bm M_c  \bm M_c^{-1} \bm P = \bm P$, which  proofs relationship \eqref{eq:MMcinv}.

Now, consider the characteristic equation of the
constraint mass matrix
\begin{equation} \notag
\big( \bm P \bm M \bm P + \nu(\bm I - \bm P) \big) \bm x - \lambda
\bm x =0
\end{equation}
Clearly $\lambda=\nu$ is the eigenvalue for all orthogonal
eigenvectors which span ${\cal N}^{\perp}(\bm A)$ because $\lambda=\nu$ means $(\bm P \bm M \bm P - \bm P)\bm x
=\bm 0 \quad \forall \bm x \in {\cal N}^{\perp}(\bm A)$. The
remaining set of orthogonal eigenvectors must lie in ${\cal N}(\bm
A)$ that are corresponding to the non-zero eigenvalues of $\bm P
\bm M \bm P$
\begin{equation} \notag
\bm P \bm M \bm P \bm x - \lambda \bm x =\bm 0 \qquad  \lambda
\neq 0 \quad \forall \bm x \in {\cal N}
\end{equation}
Therefore, the set of all eigenvalues of the p.d. matrix ${\bm
M}_c$ is  the union of the above sets corresponding to the
eigenvectors in ${\cal N}$ and ${\cal N}^{\perp}$, i.e.,
\begin{equation} \label{eq:lamb(M_bar)}
\lambda({\bm M}_c)=: \big\{\underbrace{\nu, \cdots, \nu}_{r},
\; \underbrace{ \lambda_{\stackrel{\rm min}{\neq0}} (\bm P \bm M
\bm P), \cdots, \lambda_{\rm max}(\bm P \bm M \bm P)}_{n-r} \big\}
\end{equation}
where $\{ \lambda_{\stackrel{\rm min}{\neq0}}(\bm P \bm M \bm P),
\cdots, \lambda_{\rm max}(\bm P \bm M \bm P) \}$ are all non-zero
eigenvalues of $\bm P \bm M \bm P$. According to
\eqref{eq:lamb(M_bar)} the condition number of $\bm M_c$,
which is simply the ratio of the largest to smallest eigenvalues,
is
\begin{equation} \label{eq:cond_M}
\mbox{cond}({\bm M}_c) = \frac{\max(\nu, \lambda_{\rm max}(\bm
P \bm M \bm P) )}{\min(\nu, \lambda_{\stackrel{\rm min}{\neq0}
}(\bm P \bm M \bm P))}
\end{equation}
Clearly, the RHS of \eqref{eq:cond_M} is at its minimum if $\nu$
is selected to be within the lower- and upper-bounds defined in \eqref{eq:condMc}.

By virtue of \eqref{eq:lamb(M_bar)}, the Singular Value Decomposition of $\bm M_c$ takes the form
\begin{equation}
\bm M_c = \begin{bmatrix} \bm V_1 & \bm V_2 \end{bmatrix} \begin{bmatrix} \nu \bm I  & \bm 0 \\ \bm 0 & \bm\Sigma \end{bmatrix} \begin{bmatrix} \bm V_1^T \\ \bm V_2^T \end{bmatrix}
\end{equation}
where matrix $\bm\Sigma=\mbox{diag}\{\lambda_{\stackrel{\rm min}{\neq 0}}(\bm P \bm M \bm P),\cdots,\lambda_{\rm max}(\bm P \bm M \bm P)\}$ contains the non-zero singular
values, $\bm V=[\bm V_1 \;\; \bm V_2]$ is a unitary matrix so that $\mbox{span}(\bm V_1) \equiv {\cal N}^{\perp}$ and $\mbox{span}(\bm V_2) \equiv {\cal N}$, i.e., $\bm P = \bm V_2 \bm V_2^T$, $\bm V_2^T \bm V_2 =\bm I$, and $\bm V_1^T \bm V_2 =\bm 0$. Thus
\begin{align} \notag
\bm M_c^{-1} \bm P &= \begin{bmatrix} \bm V_1 & \bm V_2 \end{bmatrix} \begin{bmatrix} \nu^{-1} \bm I & \bm 0 \\ \bm 0 & \bm\Sigma^{-1} \end{bmatrix} \begin{bmatrix} \bm V_1^T \\ \bm V_2^T  \end{bmatrix} \bm V_2 \bm V_2^T \\ \notag &= \bm V_2 \bm\Sigma^{-1} \bm V_2^T = \bm M_o^+
\end{align}

\subsection{Properties of $\bm S$} \label{appx:S_properties}
Using \eqref{eq:PMP} in the following derivations yields
\begin{align*}
\bm S \bm P &= \bm P - \bm M_c^{-1} \bm P \bm M \bm P\\
&= \bm P - \bm M_c^{-1} \bm M_c \bm P\\
& = \bm P - \bm P = \bm 0,
\end{align*}
which proves identity \eqref{eq:SP=0}. It follows
\begin{equation} \notag
(\bm I - \bm P) \bm S^T = \bm S^T - \bm P \bm S^T = \bm S^T
\end{equation}

On the other hand, by definition we have
\begin{align*}
(\bm I - \bm P) \bm S &= \bm I - \bm P - \bm M_c^{-1} \bm P \bm M + \bm P \bm M_c^{-1} \bm P \bm M\\
&= \bm I  - \bm P - \bm M_c^{-1} \bm P \bm M + \bm M_c^{-1} \bm P^2 \bm M\\
& = \bm I - \bm P,
\end{align*}
which proves identity \eqref{(I-P)St}.

Finally
\begin{align*}
\bm M \bm S &= \bm M - \bm M \bm M_c^{-1} \bm P \bm M\\
&= \bm M - \bm M \bm P \bm M_c^{-1} \bm M\\
&= (\bm I - \bm M \bm P \bm M_c^{-1}) \bm M\\
& = \bm S^T \bm M
\end{align*}
On the other hand, using the above result in the following derivation yields
\begin{equation*}
\bm S^T \bm M \bm S = \bm S^T \bm S^T \bm M = \bm S^T \bm M,
\end{equation*}
which proves identity \eqref{MS=StM}.

\bibliographystyle{IEEEtran}


\end{document}